\documentclass[11pt]{amsart}

\usepackage[utf8x]{inputenc}

\usepackage{amsmath, amsthm, amsfonts, amssymb}
\usepackage{epsfig, enumerate}
\usepackage{color}

\usepackage{hyperref}
\usepackage[all]{xy}

\usepackage{float}

\usepackage{tikz}
\usetikzlibrary{arrows,shapes,trees}
\usepackage{graphicx}
\usepackage{amssymb}
\usepackage{amsfonts}
\usepackage{amsthm}
\usepackage{amsmath}
\usepackage{mathrsfs}

\newtheorem{maintheorem}{Theorem} 
\newtheorem{theorem}{Theorem}
\newtheorem{lemma}[theorem]{Lemma}
\newtheorem{coro}[theorem]{Corollary}
\newtheorem{prop}[theorem]{Proposition}

\newtheorem*{assumptions*}{Assumptions}

\newtheorem*{rem*}{Remark}

\theoremstyle{remark}
\newtheorem{remark}[theorem]{Remark}
\newtheorem*{remark*}{Remark}

\theoremstyle{definition}
\newtheorem{definition}{Definition}

\newcommand{\A}{{\mathbf A}}
\newcommand{\B}{{\mathbf B}}
\newcommand{\C}{{\mathbf C}}

\newcommand{\G}{{\mathbf G}}

\newcommand{\V}{{\mathbf V}}
\newcommand{\W}{{\mathbf W}}

\newcommand{\Z}{{\mathbf Z}}

\newcommand{\Hh}{{\mathcal H}}

\newcommand{\CZ}{{\mathcal C}{\mathcal Z}}

\newcommand{\NN}{{\mathbb N}}

\newcommand{\RR}{{\mathbb R}}
\newcommand{\Sb}{{\mathbb S}}
\newcommand{\TT}{{\mathbb T}}

\newcommand{\ZZ}{{\mathbb Z}}

\newcommand{\be}[1]{\begin{equation} \label{#1} }
\newcommand{\ee}{\end{equation}}
\newcommand{\beq}{\begin{equation}}

\def \Diff{{\rm Diff}}

\def \W{{\mathcal W}}

\def \A{{\mathcal A}}

\def \dist{{\mathrm{dist}}}
\def \C{{\mathcal{C}}}

\def \W{\mathcal{W}}
\def \cW{\mathcal{W}}

\def \dist{{\mathrm{dist}}}
\def \Ho{{\mathrm{Homeo}}}

\def \hx0{\hat{x_0}}

\def \id{\mathrm{id}}

\def \G{\mathcal{G}}
\def \B{\mathcal{B}}
\def \V{\mathcal{V}}
\def \Z{\mathcal{Z}}

\begin{document}

\title[Transitive Centralizer]{Transitive Centralizer and fibered partially hyperbolic systems}

\author{Danijela Damjanovi\'c}

\address[Damjanovi\'c]{Department of mathematics, Kungliga Tekniska högskolan, Lindstedtsv\"agen 25, SE-100 44 Stockholm, Sweden.} 
\email{ddam@kth.se}

\author{Amie Wilkinson}
\address[Wilkinson]{Department of mathematics, the University of Chicago, Chicago, IL, US, 60637}
\email{wilkinso@math.uchicago.edu}

\author{Disheng Xu}
\address[Xu]{Great bay University, Songshanhu International Community,  Dongguan, Guangdong, China, 523000}
\email{xudisheng@gbu.edu.cn}

\maketitle
\begin{abstract}
   We prove several rigidity results about the centralizer of a smooth diffeomorphism, concentrating on two families of examples: diffeomorphisms with transitive centralizer, and perturbations of isometric extensions of Anosov diffeomorphisms of nilmanifolds. 

   We classify all smooth diffeomorphisms with transitive centralizer: they are exactly the maps that preserve a principal fiber bundle structure, acting minimally on the fibers and trivially on the base. 

   We also show that for any smooth, accessible isometric extension 
   $f_0\colon M\to M$ of an Anosov diffeomorphism of a nilmanifold, subject to a spectral bunching condition,  any  $f\in \Diff^\infty(M)$ sufficiently $C^1$-close to $f_0$ has centralizer a Lie group. If the dimension of this Lie group equals the dimension of the fiber, then  $f$ is a principal fiber bundle morphism covering an Anosov diffeomorphism.

   Using the results of this paper, we further classify the centralizer of any partially hyperbolic diffeomorphism on a $3$-dimensional, nontoral nilmanifold: either the centralizer is virtually trivial, or the diffeomorphism is an isometric extension of an Anosov diffeomorphism, and the centralizer is virtually $\ZZ\times \TT$.
\end{abstract}

\tableofcontents

\section{Introduction}\label{intro}

Let $M$ be a closed, connected smooth manifold. The {\em centralizer} $\Z(f)$ of a diffeomorphism $f\colon M\to M$ is the set of all  diffeomorphisms that commute with $f$ under composition.  The centralizer may be regarded as the set of smooth symmetries of $f$. The group $\Z(f)$ always contains the iterates of $f$ as a normal subgroup: $\Z(f) \rhd \langle f \rangle$.  For the $C^1$-generic  diffeomorphism $f$ of a closed manifold, the centralizer is no bigger than its iterates, i.e. $\Z(f) = \langle f \rangle$ \cite{BCW}.  In the latter case, we say that $f$ has {\em trivial centralizer}. Whether this genericity of trivial centralizers remains true in higher regularity classes (say $C^2$)  remains an open question, except on the circle \cite{Kopell}.

In this paper and our earlier work \cite{DWXcent} we address the general question: 

\medskip

\centerline{\em If $\Z(f)$ is nontrivial, what can we say about $f$?}

\medskip

Consider the most extreme situation in which the centralizer $f$ is $\Diff(M)$.  In this case, if there is a point $x$ that is not fixed by $f$, then there is an element $g\in \Z(f)$ fixing $x$ and {\em not} fixing $f(x)$. This leads to a contradiction, since $g(f(x)) = f(g(x)) = f(x) = x$.  Thus $\Z(f) = \Diff(M)$ implies that $f=
Id_M$.  Indeed, as long as $\Z(f)$ acts doubly transitively on $M$, the same conclusion holds.\footnote{All this argument requires is that for every $x,y\in M$, there exists $g\in \Z(f)$ such that $g(x)=x$ and $g(y)\neq y$.}

Even the weaker case in which $\Z(f)$ acts {\em transitively} on $M$ puts strong constraints on $f$.  For example if $f$ is minimal and $\Z(f)$ acts transitively, then $f$ is a minimal translation on a torus (after a smooth conjugacy).  More generally, we show:\footnote{For $r\in (0,1)$, it is possible there is a similar result assuming only H\"older regularity. For $r=0$, this appears to be open at least in dimension at least 2.}

\begin{maintheorem}\label{tmain=transitivecent}
Let $f$ be  homeomorphism of a connected closed manifold $M$, and denote by $\Z^r(f)$ the $C^r$ centralizer of $f$ for $r\in \ZZ_{\geq 1}$. Then $\Z^r(f)$ acts transitively on $M$ if and only if $M$ is a $C^r$ principal $\mathbb T^k\times \ZZ/d\ZZ$ fiber bundle for some $k\geq 0, d\geq 1$, and $f$ acts as the identity on the base and a constant minimal translation on the fibers. 

Moreover if $\Z^r(f)$ acts transitively on $M$, then $f$ is minimal on $M$ if and only if  $f$ is $C^r$ conjugate to a minimal translation on a torus.
\end{maintheorem}

Thus {\em a sufficiently rich centralizer characterizes some systems of low complexity.}

\medskip

At the other extreme of dynamical complexity are the Anosov diffeomorphisms. For an irreducible Anosov diffeomorphism $f$ of a torus\footnote{For a toral Anosov diffeomorphism $f$, {\em irreducibility} means that the linearization $f_0$ of $f$ does not preserve a proper invariant subtorus.}, there are only finitely many possibilities for the centralizer (up to isomorphism and finite index subgroups), all free abelian. More generally, the centralizer of a transitive Anosov diffeomorphism is always discrete.

In this paper, we also consider systems that combine features of low and high complexity, namely the isometric extensions of Anosov diffeomorphisms of nilmanifolds (and their perturbations). Such extensions are given by a smooth fiber bundle $F\hookrightarrow M \to N$, where $M$ is closed and connected and $N$ is a nilmanifold,
and a diffeomorphism $f_0\in \hbox{Diff}^\infty(M)$ preserving the fiber bundle structure, acting isometrically on fibers, and projecting to an Anosov diffeomorphism $\bar f_0$ of $N$, see \cite{Bo}.

We assume a generic condition on $f_0$ called {\em accessibility} and that $\bar f_0$ satisfies a spectral condition (satisfied, for example, by an affine diffeomorphism of $M$). For precise definitions and  further discussion, see Section \ref{sec: pre}.   Under these assumptions, we
obtain a type of centralizer rigidity, whose proof uses elements of the proof of Theorem \ref{tmain=transitivecent}.




\begin{maintheorem}\label{main: application iso ext}Let $f_0\in \Diff^\infty(M)$ be an accessible isometric extension of an Anosov diffeomorphism $\bar{f}_0\colon N\to N$ of a closed nilmanifold $N$. Assume that $D\bar{f_0}|_{E^s_{\bar{f_0}}}$ and $D\bar{f_0}^{-1}|_{E^u_{\bar{f_0}}}$ have narrow spectrum. Then there exists $s_0 = s_0(f_0)>0$ such that for any $s\geq s_0$ and any $f\in \Diff^\infty(M)$ that is  sufficiently $C^1$-close to $f_0$,
the group $\Z^s(f)$ is a $C^0$-closed, $k$-dimensional Lie subgroup of  $\Diff^s(M)$,  with $k\leq \dim 
 E^c_{f_0}$.  Moreover, at least one of the following possibilities holds:

\begin{enumerate}
\item  $k=0$, and $\Z^s(f)$ is a countable discrete subgroup of $\Ho(M)$, or
\item $k>0$, and there exist a $k$-dimensional compact Lie group $G$ and $r\gg s$  such that $M$ admits a $C^r$, $f$-invariant principal $G$-bundle structure that is subordinate to $\W^c_f$. Along each $G$-fiber $f$ acts as left translation, and $\Z^s(f)$ is a countable extension by the group ($\cong G$) formed by all the right translations along $G$-fibers.  
\item  $k=\dim E^c_{f_0}$, and $f$ is an automorphism of a $C^\infty$ principal fiber bundle projecting to an Anosov diffeomorphism with narrow   spectrum. 
\end{enumerate}
\end{maintheorem}
\begin{remark}Theorem \ref{main: application iso ext} also holds with  ``narrow  spectrum" replaced by  ``(pointwise) $\frac{1}{2}$-pinched" -- see Section~\ref{ss=normalforms} for definitions.
\end{remark}
\begin{remark}
    Note that in Theorem~\ref{main: application iso ext} (in contrast to Theorem~\ref{tmain=transitivecent}),  the group $G$ might not be abelian, as the following example shows. Fix a compact group  $G$, and an Anosov diffeomorphism $\bar f\colon N\to N$ of a closed manifold $N$. Let $\theta\colon N\to G$ be a smooth function. The manifold $M=N\times G$ is a principal $G$-bundle with respect to the right action $g\cdot (x,h) = (x,hg)$. The diffeomorphism $\bar f_\theta\colon M\to M$ defined by $\bar f_\theta(x,h) = (\bar f(x), \theta(x) h)$ is partially hyperbolic and commutes with the right action of $G$. 
    
    Theorem B of \cite{BW99} shows that if $N$ is nilmanifold, then for the generic such $\theta$ (in a very strong sense), the map $\bar f_\theta$ is also accessible.  Indeed \cite[Theorem B]{BW99} classifies all non-accessible examples: for example, if $G$ is semisimple, and $\bar f_\theta$ is ergodic,  then $\bar f_\theta$ is  accessible.
\end{remark}
 \begin{remark}  Since any infinite compact Lie group contains subgroups isomorphic to $\TT^1$, Theorem \ref{main: application iso ext} has the following consequence: if $\W^c_f$ does not admit a  free $\TT^1$ action, then
$\Z^s(f)$ is $C^0$-discrete subgroup of  $\Diff^s(M)$.  For example, any accessible isometric extension of an Anosov  automorphism  with $\Sb^2$ fiber has discrete countable centralizer. More generally, if the fiber has non-zero Euler characteristic, then the centralizer is discrete. 
\end{remark}
\begin{remark}
The centralizer of a diffeomorphism in general is not $C^0$-closed.  For example, consider a diffeomorphism $f$ of the circle with Liouville rotation number, which is not $C^1$ conjugate to a rotation.  Then $\Z^2(f)=\langle f \rangle$, which is not closed in $\Ho(\Sb^1)$ but is discrete, and hence closed, in $\Diff(\Sb^1)$.
\end{remark}
Theorem \ref{main: application iso ext} and its variants have various applications in the study of the centralizer of partially hyperbolic systems. For one example, it is now possible to strengthen some of the results in \cite{DWXcent} to remove the volume-preserving hypotheses.
In this prior work, we established centralizer rigidity for some volume-preserving partially hyperbolic systems with a $1$-dimensional center foliation, including certain discretized Anosov flows (and Barthelm\'e and Gogolev \cite{BG} extended this work to discretized Anosov flows on $3$-manifolds). 

Recently, Wendy Zhijing Wang \cite{ZW} has generalized Theorem \ref{main: application iso ext} to the setting of center-fixing partially hyperbolic systems, and her work can be used to remove the volume preserving assumptions for the results on discretized Anosov flows in \cite{DWXcent, BG}.  In particular, \cite{ZW} obtains as a corollary, using \cite{BG}, that for any $3$-manifold $M$ with nonsolvable fundamental group supporting an Anosov flow $\varphi_t$, any $C^1$-small, smooth perturbation of $\varphi_1$ either has virtually trivial centralizer or embeds in a smooth Anosov flow.

In a similar vein, combining techniques in this paper with work of Hammerlindl and Potrie \cite{HP}, we obtain a global characterization of  the centralizer of a partially hyperbolic diffeomorphism of any nontoral, $3$-dimensional nilmanifold.
\begin{maintheorem}\label{coro: global Heis3}Let $f\in \Diff^\infty(M)$ be a  partially hyperbolic diffeomorphism of a $3$-dimensional nilmanifold $M$, $M\neq \TT^3$. Then $\Z^\infty(f)$ is either virtually trivial or virtually $\ZZ\times \TT^1$. In the latter case, $f$ is a $C^\infty$ isometric extension of a toral Anosov diffeomorphism.
\end{maintheorem}

\paragraph{\textbf{Acknowledgements:}}D.D. is supported by Swedish Research Council grant 2019-67250. 
 A.W. is supported by NSF grants DMS-1900411 and DMS-2154796.
D.X. is supported by NSFC grant 12090015. 

Part of this work was done during the third author's visit to Kungliga Tekniska Högskolan and the University of Chicago. D.X. thanks Kungliga Tekniska Högskolan and the University of
Chicago for hospitality during  his visits. 

The authors thank Artur Avila for many useful discussions during an early stage of this project.
\section{Preliminaries}\label{sec: pre}
\subsection{Partially hyperbolic dynamics} We recall some definitions from smooth dynamics. Let $M$ be a complete Riemannian manifold, and let $f\in \Diff(M)$.  A {\em dominated splitting} for $f$ is a
direct sum decomposition of the tangent bundle
\[TM=E_f^1\oplus E_f^2\oplus\cdots \oplus E_f^k\]
such that 
\begin{itemize}
\item the bundles $E_f^i$ are {\em $Df$-invariant}: for every $i\in \{1,\ldots, k\}$ and $x\in M$, we have $D_xf(E_f^i(x)) = E_f^i(f(x))$; and
\item  $Df\vert_{E_f^i}$ {\em  dominates} $Df\vert_{E_f^{i+1}}$: there exists $N\geq 1$ such that for any $x\in M$ and any unit vectors $u\in E_f^{i+1}$ , and $v\in  E_f^{i}$:
\[ \|D_xf^N (u)\| \leq  \frac{1}{2} \,\|D_x f^N(v)\|.\]
\end{itemize}

A $C^{1}$ diffeomorphism $f: M \rightarrow M$ of a complete Riemannian manifold $M$ is \emph{partially hyperbolic} if there is a dominated splitting $TM = E_f^{u} \oplus E_f^{c} \oplus E_f^{s}$ and $N \geq 1$ such that for any $x \in M$, and any choice of unit vectors $v^{s} \in E_f^{s}(x)$ and $v^{u} \in E_f^{u}(x)$, we have
\[
\max\{\|D_xf^{N}(v^{s})\|,  \|D_xf^{-N}(v^{u})\|\} < 1/2.
\]
We always assume the bundles $E_f^{s}$ and $E_f^{u}$ are nontrivial. If $E_f^c$ is trivial then $f$ is  \emph{Anosov}.  


In many cases of interest here, we consider a partially hyperbolic diffeomorphism $f$ equipped with a {\em center foliation $\W_f^c$}  that is tangent to $E_f^c$ and whose leaves are compact, forming a fibration.  We distinguish between several cases of such fibered systems.

\begin{definition}\label{def: fiber systems}Let $f$ be a partially hyperbolic diffeomorphism of a closed manifold $M$. Assume that there exists an $f-$invariant center foliation $\W^c_f$ with compact leaves.
\begin{itemize}
\item If $\W^c_f$ is a topological fibration of $M$, i.e. the quotient space $M/\W^c_f$ is a topological manifold\footnote{Or, equivalently, if $\W^c_f$ has trivial holonomy; see \cite{Bo}}, then $f$ is called a \emph{fibered partially hyperbolic system}, and the map $\bar{f}:M/\W_f^c\to M/\W_f^c$ canonically induced by $f$ is called the \emph{base map}.
\item A fibered partially hyperbolic system $f$ is  \emph{smoothly fibered} (or \emph{$C^{r}$ fibered}, for $r\geq 1$ )  if $\W^c_f$ is a $C^\infty$ (respectively $C^r$) foliation, and $f$ is $C^\infty$ (resp. $C^r$). In this case, the base map is an Anosov diffeomorphism.
\item  A smoothly fibered partially hyperbolic system $f$ is an \emph{isometric extension of an Anosov diffeomorphism}, if there is a $C^\infty$, $f-$invariant metric on $E^c_f$.

\item A smoothly fibered partially hyperbolic system $f$ is a \emph{compact Lie group extension of an Anosov diffeomorphism} if there is a compact Lie group $G$ and a $C^\infty$ $f-$invariant principle $G$-bundle structure $\pi:M\to M/\W^c_f$ such that, restricted to every fiber, $f$ acts as a left $G$-translation.
\end{itemize}
\end{definition}

\begin{remark}For the last two items in Definition \ref{def: fiber systems}, the partial hyperbolicity assumptions are not necessary, since any isometric extension of an Anosov diffeomorphism is partially hyperbolic; see \cite[Proposition 5.3]{Dou}.
\end{remark}

If $f$ is partially hyperbolic,  then there are foliations $\W^u_f$
and $W^s_f$ tangent to the bundles $E^u_f$ and $E^c_f$, respectively.
\begin{definition} We say that $x, y\in M$ are in the
same {\em accessibility class} if they can be joined by an {\em $su$-path}, that is, a piecewise $C^1$ path such that every arc (or {\em leg}) is contained in a single leaf of $\W^s_f$ or a single leaf of $\W^u_f$. We say that $f$ is {\em accessible} if $M$ consists of a single accessibility class. 
\end{definition}

\begin{definition}
For $r \geq 1$, we say that a partially hyperbolic diffeomorphism $f$ of a Riemannian manifold $M$ is   {\em $r-$bunched}   if there exists $k\geq  1$  such that:
\[\sup_p \left \{ \|D_p f^k|_{E^s}\|\cdot \|(D_p f^k|_{E^c})^{-1}\|^r, \|(D_p f^k|_{E^u})^{-1}\|\cdot \|D_p f^k|_{E^c}\|^r\right \}<1,\]
\begin{eqnarray*}
&&\sup_p \|D_p f^k|_{E^s}\|\cdot \|(D_p f^k|_{E^c})^{-1}\|\cdot \|D_p f^k|_{E^c}\|^r<1, \hbox{ and}\\
&&\sup_p \|(D_p f^k|_{E^u})^{-1}\|\cdot \|D_p f^k|_{E^c}\|\cdot \|(D_p f^k|_{E^c})^{-1}\|^r<1.
\end{eqnarray*}
If $f$ is an isometric extension of an Anosov diffeomorphism, then $f$ is $\infty$-bunched (i.e. $r$-bunched, for all $r$).  For any such $f$, and any finite $r>1$, if $\hat f$ is $C^r$ and sufficiently $C^1$-close to $f$, then $\hat f$ is $r$-bunched.

If $f$ is an $r$-bunched, $C^{r}$ fibered partially hyperbolic diffeomorphism, then the leaves of  $\W_f^c, \W_f^{cs},\W_f^{cu}$ are $C^r$, moreover, if $f$ is $C^{r+1}$, then the stable and unstable \emph{holonomies} and the bundles $E_f^s, E_f^u$ are $C^r$ along $\W_f^c$; see \cite{PSW, W}.

\end{definition}

\begin{definition}
Denote by $\CZ^s(f)$ the group of all $g\in \Z^s(f),s=0,1,\dots$ fixing the center leaves of $f$; that is, $g\in \CZ^s(f)$  if and only if $g\W^{c}(x)=\W^c(x)$, for all $x\in M.$ \end{definition}

\subsection{$su$-holonomy}\label{subsec: suholon}

Every  fibered partially system $f\in \Diff(M)$ admits {\em global $su$-holonomies}, meaning that for any $su$-path 
$\gamma$ in $M$ from $x$ to $y$, there exists a unique homoemorphism $H_\gamma\colon \W^c_f(x)\to \W^c_f(y)$ with the properties:
\begin{itemize}
\item $H_\gamma(x) = y$;
\item $H_{\gamma_1\cdot \gamma_2} = H_{\gamma_2}\circ H_{\gamma_1}$, where $\gamma_1\cdot \gamma_2$ is the concatenation of two $su$-paths; and
\item if $\gamma$ is tangent to $\W^s_f$ (resp $\W^u_f$), then $H_\gamma$ is the $\W^s_f$ (resp $\W^u_f$) holonomy, restricted to $\W^c_f$ leaves.
\end{itemize}


The set of $su$-holonomies from a fixed center leaf $\W^c_f(x_0)$ to itself form a group, which we denote by $\mathcal{H}_{f}(x_0)$.  If $f$ is $r$-bunched, then $\mathcal{H}_{f}(x_0)<\Diff^{r}(\W^c_f(x_0))$.

\subsection{Leafwise structural stability}\label{subsec: fbratn}


In the setting of fibered partially hyperbolic systems, there is a variety of results we will use, starting with the Hirsch-Pugh-Shub perturbation theory  \cite[Theorems 7.5 and 7.6]{HPS} (see Remark 4 on p. 117), \cite[Theorem 7.1]{HPS}, and \cite[Theorems A and B]{PSW}. See also the discussion in \cite[Section 3]{PSW}. 
\subsection{Stable accessibility}In this section we establish the stable accessibility of any isometric extension of an Anosov diffeomorphism. 
\begin{prop}\label{prop: stab acc}Let $f\in \Diff^\infty(M)$ be an  isometric extension of an Anosov  diffeomorphism. If $f$ is accessible, then $f$ is stably accessible, i.e.~any $g$ sufficiently $C^1$-close to $f$ is accessible.
\end{prop}
\begin{proof}We follow the strategy of \cite{BW99}. Since most of the proofs are the same or similar, we only sketch the proof here and recommend reading \cite{BW99} for background.  Fix a center leaf $\W^c_f(x_0)$.
\\
\\
\textbf{Step 1: accessibility implies local accessibility.}  
Any $f$-invariant $G$-structure on the leaves of $\cW_f^c$ is also invariant under the stable and unstable holonomies between $\W_f^c$ leaves;   since $f$ is an isometric extension,  all the $su$-holonomies from $\W^c_f(x_0)$ to itself are isometries of $\W^c_f(x_0)$. Let  $\Hh = \Hh_{f}(x_0)$ be the group all such holonomies, and denote by $\Hh^0$ the subgroup of $\Hh$ consisting of all the holonomies $H_\gamma$ with $\pi\circ\gamma$ a closed, contractible path in the base. Note that such a $\pi\circ\gamma$ is also an $su$-path,  for the base Anosov diffeomorphism $\bar f$ and that $H_\gamma$ depends only on the projection $\pi\circ\gamma$.

Then $\Hh$ is a subgroup of the Lie group $\mathrm{Iso}(\W^c_f(x_0))$, and $\Hh^0$ is the identity component of $\Hh$. Any element $H_\gamma$ of $\Hh^0$ can be isotoped to the identity through maps $H_{\gamma_t}$ with $\pi\circ\gamma_t$ an isotopy to the trivial path through closed $su$-paths for $\bar f$ (``Brin's argument"). Thus $\Hh^0$ is path connected. A result of Kuranashi-Yamabe \cite{Y} then implies that $\Hh^0$ is  a Lie subgroup of  $\mathrm{Iso}(\W^c_f(x_0))$.

The quotient $\Hh/\Hh^0$ has at most countably many elements (since each coset corresponds to an element of $\pi_1(M/\W^c_f, \bar{x_0})$). The Baire category theorem then implies that $\mathrm{Image}(\Hh^0\cdot x_0)$ has nonempty interior. Since $\Hh^0$ is a Lie subgroup of the isometry group of $\W_f^c(x_0)$, this implies that $\mathrm{Image}(\Hh^0\cdot x_0) = \W_f^c(x_0)$. Hence $f$ is accessible through those $su$-paths that project to null-homotopic, closed paths in $M/\W^c$.
\\
\\
\textbf{Step 2: a homogeneous space structure on the leaves of the center foliation.}  From Step 1, we have that $\Hh^0$ is a connected Lie group acting transitively on $\W^c_f(x_0)$. This extends to a transitive action of $\Hh^0$ on every center leaf via conjugation (for each $x\in M$, one fixes an $su$-holonomy between $\W^c_f(x)$ and  $\W^c_f(x_0)$ and conjugates the $\Hh_0$-action by this map).  Thus $\cW^c_f$ is an $\Hh^0$-bundle; i.e., there is an $f$-invariant homogeneous space structure on the foliation $\W^c_f$ obtained by identifying each $\W^c_f(x)$ with $\Hh^0(x)/\mathrm{Stab}(x)$, where $\mathrm{Stab}(x)$ is the stabilizer of $x$.
\\
\\
\textbf{Step 3: a useful criterion for stable accessibility.} The paper \cite{BW99} contains a useful criterion for stable accessibility of dynamically coherent partially hyperbolic diffeomorphisms, which will be used in our proof.   We say that a partially hyperbolic diffeomorphism $f$ is {\em dynamically coherent} if there exist $f$-invariant center stable and center unstable foliations $\W^{cu}_f$ and $\W^{cs}_f$ tangent to the bundles $E_f^c\oplus E_f^u$ and  $E_f^c\oplus E_f^s$, respectively. The fibered partially hyperbolic maps we consider here are dynamically coherent (see Theorem 8 in \cite{DWXcent}).

Let $f$ be a dynamically coherent partially hyperbolic diffeomorphism, and let $d_c$ be the dimension of $E^c_f$. A point $q_0\in M$ is \textit{centrally engulfed} from a point $p_0\in M$ if there is a continuous map $\Gamma: Z \times [0, 1] \to M$ such that:
\begin{enumerate}
    \item $Z$ is a compact, connected, orientable, $d_c$-dimensional manifold with boundary;
    \item for each $z\in Z$, the curve $\gamma_z(\cdot ) = \Gamma(z,\cdot)$ is an $su$-path with $\gamma_z(0) = p_0$ and $\gamma_z(1)\in \W^c_f(q_0)$;
    \item there is a constant $C$ such that every path $\gamma_z$ has at most $C$ legs;
    \item $\gamma_z(1)\neq q_0, \forall z\in \partial Z$;
    \item the map $(Z, \partial Z) \to  (\W^c_f(q_0), \W^c_f(q_0)-\{q_0\})$ defined by $z\to  \Gamma (z, 1)$ has positive degree.
\end{enumerate}
We have the following criterion for stable accessibility.
\begin{theorem}[Corollary 5.3 in \cite{BW99}] \label{thm: crit stab acc}
Let $f$ be a dynamically coherent partially hyperbolic
diffeomorphism. Suppose that $f$ is accessible and that exist $p_0\in M$ such that $p_0$ can be centrally engulfed from $p_0$. Then $f$ is stably accessible.
\end{theorem}
\paragraph{\textbf{Step 4: achievable and approximable paths.}} By Theorem \ref{thm: crit stab acc}, it suffices to show $x_0$ is centrally engulfed from $x_0$. We follow the strategy of Section 9 of \cite{BW99}: the key is to use the homogeneous space structure of $\W^c_f(x_0)$ to construct the map $\Gamma$. We say that  $\psi\colon Z\to \Hh^0$ is \textit{achievable} if it is the “endpoint map” of a continuous family of $su$-holonomies that begin and end in $\W^c_f(x_0)$ and project to null-homotopic loops in $M/\W^c_f$ (i.e. through elements of $\Hh^0_{x_0}(f)$).

More precisely, $\psi$ is  achievable if there exist a continuous function \[\Xi=(\Xi_1,\Xi_2) \colon Z\times [0, 1] \to \{(\bar{x},g_x):\bar{x}\in M/\W^c_f, g_x\in \mathrm{Iso}(\W^c_f(x_0), \W^c_f(x))\}\] and an integer $C\geq 1$ such that 
\begin{enumerate}
    \item for each $z\in Z$, $\Xi_1(z,t)$ is an $su$-path on the base $M/\W_f^c$  beginning and ending in $\bar{x_0}$, with at most $C$ legs;
    \item for each $z\in Z$ and $t\in [0,1]$, 
    $\Xi_2(z,t)$ is the (isometric) $su$-holonomy induced by lifting the $su$-path  $\Xi_1(z,t)$; and
    \item the path $\xi_z :[0,1]\to M$ defined by $\xi_z(t)=\Xi_2(z,t)$ has the property that $\xi_z(1)$ maps any $y\in \W^c_f(x_0)$ to $\psi(z)\cdot y$.
    \end{enumerate}


Following \cite{BW99}, we say a map $\psi\colon Z\to \Hh^0$ is \textit{approximable} if for each $\epsilon>0$ there is an achievable map $\psi^\epsilon : Z \to  \Hh^0$ such that $\dist_{C^0} (\psi, \psi^\epsilon) < \epsilon$. 
\begin{lemma}[Lemma 9.2 of \cite{BW99}]\label{lemma: prod apprx}
Let $\psi_i \colon Z_i \to G, 1\leq i \leq k$ be approximable maps. Then the product map
\[(z_1,\dots,z_k) \mapsto  \psi_k(z_k)\psi_{k-1}(z_{k-1})\cdots \psi_2(z_2)\psi_1(z_1)\]
is approximable.    
\end{lemma}
\paragraph{\textbf{Step 5: geodesics in $\Hh^0$ are approximable.}} By the same proof as Proposition 9.3 of \cite{BW99}, we have that for any $v$ in the Lie algebra of $\Hh^0$, the map $\sigma_v\colon [-1,1]\to \Hh^0, z\mapsto \exp(zv)$ is approximable. Hence by Lemma \ref{lemma: prod apprx}, for any $v_1,\dots, v_{d_c}$ in the Lie algebra of $\Hh^0$, the map $\psi_{v_1,\dots, v_{d_c}}\colon [-1,1]^{d_c}\to \Hh^0$ defined by $(z_1,\dots, z_{d_c})\mapsto \prod^1_{i=d_c} \exp(z_iv_i)$ is approximable. Let $\langle\cdot,\cdot\rangle$ be the invariant inner product on the Lie algebra of $\mathrm{Iso}(\W^c_f(x_0))$. With respect to   $\langle\cdot,\cdot\rangle$ , we choose the $v_i$ such that 
\begin{enumerate}
    \item $\langle v_i,v_j\rangle=\eta\delta_{ij}$, for some small $\eta$, so that $\exp$ is ``almost an isometry" on the $\eta$-neighborhood of $0$; and
    \item the angle between the subspace spanned by the $v_i$ and the Lie algebra of $\mathrm{Stab}(x_0)$ is greater than $1/10$.
\end{enumerate}
Then $\psi_{v_1,\dots,v_{d_c}}$ is a  parametrization of a $d^c$-dimensional submanifold of $\Hh^0$, whose interior contains the identity element $id$ and is uniformly transverse to $\mathrm{Stab}(x_0)$. Since $\psi_{v_1,\dots,v_{d_c}}$ is approximable, it has an achievable approximation $\psi\colon[-1,1]^{d_c}\to \Hh^0$.  Let $Z= [-1,1]^{d_c}$, and let \[\Xi=(\Xi_1, \Xi_2)\colon Z\times [0,1]\to \{(\bar{x},g_x)\colon \bar{x}\in M/\W^c_f, g_x\in \mathrm{Iso}(\W^c_f(x_0), \W^c_f(x))\}\] be the map satisfying properties (1)-(3) in the previous step for the achievable map $\psi$, for some $C\geq 1$.

Now consider the evaluation map $\Gamma\colon [-1,1]^{d_c}\times [0,1]\to M$ defined by $\Gamma(z_1,\dots, z_{d_c},t)=\Xi_2(z_1,\dots, z_{d_c},t)(x_0)$. By the $1/10$-transversality  of the subspaces $\mathrm{Span}(v_1,\dots,v_{d_c})$ and the Lie algebra of $\mathrm{Stab}(x_0)$, if $\eta$ is sufficiently small and the achievable approximation $\psi$ is chosen sufficiently close, the map $\Gamma$ will satisfy conditions (1)-(5) in Step 3 (with $p_0, q_0$ equal to $x_0$), and so $x_0$ is centrally engulfed from $x_0$. Theorem~\ref{thm: crit stab acc} then implies that $f$ is stably accessible, completing the proof of Proposition \ref{prop: stab acc}. 
\end{proof}

\subsection{$\Z^1(f)$-invariance of $f$-invariant foliations}
If a diffeomorphism $g$ commutes with a partially hyperbolic diffeomorphism $f$, the derivative $Dg$ preserves the $Df$-invariant partially hyperbolic splitting for $f$ and the foliations $\cW^{s/u}_f$. Perhaps surprisingly, it is unknown in general whether if, under the additional assumption that $f$ is dynamically coherent, such a $g$ must also preserve even one of the foliations $\W^{cu}_f, \W^{cs}_f$, or $\W^c_f$.  Even in the case of a fibered partially hyperbolic diffeomorphism, this appears to be an open question. We do have the following result in the case that the base is a nilmanifold (and this alone accounts for our assumption that the base be a nilmanifold in Theorem~\ref{main: application iso ext}).
\begin{prop} \label{prop: fib systm zf inv} Let $f\in \Diff(M)$ be a fibered partially hyperbolic diffeomorphism such that $M/\W^c_f$ is homeomorphic to a nilmanifold $N$. Then for any $g\in \Z^{1}(f)$, we have $g\W^c_f=\W^c_f$.\end{prop}
\begin{proof}
For any fibered partially hyperbolic diffeomorphism $f$,  \cite{Bo} shows that the quotient dynamics $\bar{f}:M/\W^c_f\to M/\W^c_f$ is expansive and has the pseudo orbit tracing property (unique shadowing of pseudo-orbits).   Theorem D of \cite{Dou} then implies that $\bar{f}$ is topologically conjugate to an Anosov automorphism of $N$. From now on we identify $M/\W^c_f$ with $N$ through the conjugacy map so that $f:M\to M$ can be viewed as a bundle automorphism of a continuous fiber bundle $M \to N$ over an Anosov automorphism $\bar{f}$.

In particular, $\bar f$ has a fixed point $\bar x_0$, and so $f$ has an invariant fiber $\W^c_f(x_0) = \pi^{-1}(\bar x_0)$.
The universal cover of $N$  is a  nilpotent Lie group $G_N$.  We may assume that $\bar x_0$ lifts to the identity $0\in G_N$.  Consider the canonical covering map $p: (G_N, 0)\to (M/\W^c_f, \W^c_f(x_0))$.  The pullback fiber bundle $p^\ast(\W^c_f)$ is a covering space over the contractible space $G_N$ and hence is a trivial bundle. In particular, for every $g\in  \Z_{\Diff(M)}(f)$ there exist lifts $\hat f, \hat g$ of $f,g$ to  $p^\ast(\W^c_f)$ such that $\hat f$ commutes with $\hat g$.

To establish Proposition \ref{prop: fib systm zf inv} it suffices to  show that $\hat g$ preserves the fibration $\widehat \W^c_f$ of the bundle $p^\ast ( \W^c_f ) $. By uniform compactness of the leaves of $\W^c_f$ and $f$-invariance of $\W^c_f$, it follows that for any two points $\hat x, \hat y\in p^\ast(\W^c_f)$, we have $\hat y\in \hat \W^c_f(\hat x)$ if and only if for any $n\in \ZZ$, there is a $C^1$-path $\gamma_n$ connecting $\hat f^n(\hat x)$ to $\hat f^n(\hat y)$ whose length is bounded by some $C>0$ independent of $n$; in fact we can take $C=\max_{x\in M} \mathrm{diam}(\W^c_f(x))$. Here we use the fact the induced action of $\bar{f}$ is an Anosov  automorphism of $N$.

Now consider the points $\hat g(x), \hat g(y)$.  For any $n\in \ZZ$ the points $\hat f^n(\hat g(x))$ and $ \hat f^n(\hat g(y))$ can be linked by the path $\hat g(\gamma_n)$, where $
\gamma_n$ is defined in the previous paragraph. Moreover we have \[\text{the length of }\hat g(\gamma_n)\leq \|g\|_{C^1}\cdot \text{the length of }\gamma_n\leq \|g\|_{C^1}\cdot C,\]
where $C$ is defined in the previous paragraph and is independent of the choice of $n$. This implies that $\hat g(x), \hat g(y)$ are contained in the same $\widehat \W^c_f$ leaf; hence $\hat g$ preserves $\widehat \W^c_f$.
\end{proof}

\subsection{Normal form theory}\label{ss=normalforms} In this section we introduce some useful aspects of the normal form theory for contracting foliations from \cite{Kal-normalform}, omitting some technical details, for which we refer the reader to the source.
\begin{definition}\label{def: pinching or narrow}Let $f$ be a $C^\infty$ diffeomorphism of a compact manifold $X$ and let $\W$ be an $f$-invariant
continuous foliation of $X$ with uniformly $C^\infty$ leaves. Suppose that $\|Df|_{T\W}\| < 1$. We say that $Df|_{T\W}$ 
  \begin{enumerate}
 \item satisfies \textit{the (pointwise) $\frac{1}{2}$-pinching condition} if there exist $C > 0$ and $\gamma<1$ such that
\begin{equation*}
\|  (Df^n|_{T_x\W})^{-1}\|\cdot \|Df^n|_{T_x\W}\|^2\leq C\cdot \gamma^n, \text{ for all }x\in X, n\in \NN.
\end{equation*} 
\item \textit{has $(\chi, \epsilon)$-spectrum} if there exists an $f$-invariant dominated splitting of $T\W=\bigoplus E_i$ and a continuous Riemannian metric on $\|\cdot \|$ on $T\W$ such that 
\begin{equation*}
e^{\chi_i-\epsilon}\|v\|\leq \|Df(v)\|\leq  e^{\chi_i+\epsilon}\|v\|   \text{ for all }v\in E_i.
\end{equation*}
Here $\chi=(\chi_1,\chi_2,\dots,\chi_l)$ satisfies $\chi_1<\dots<\chi_l<0$, and we assume $\epsilon>0$.
\item \textit{has narrow  spectrum} if $Df|_{T\W}$ has $(\chi,\epsilon)$-spectrum for some $\chi_1<\dots<\chi_l<0$ and $\epsilon\in(0,\epsilon_0(\chi))$, where $\epsilon_0(\chi)$ is defined in \cite[Section 3.4]{Kal-normalform} to guarantee that the ``bands" $(\chi_i-\epsilon,\chi_i+\epsilon)$ are  narrow (see \cite{Kal-normalform} for more details).
 \end{enumerate}   
\end{definition}

The following lemma summarizes results from Section 3.4 and Theorem 4.6 of \cite{Kal-normalform}.
\begin{lemma}\label{lemma: normal form}
Let $f, X, \W$ be as in Definition \ref{def: pinching or narrow}, and suppose that  $Df|_{T\W}$ satisfies the $\frac{1}{2}$-pinching condition or has narrow spectrum.  Then there exist $d, s_0=s_0(f)\geq 1$ and a family $\{\Psi_x\}_{x\in M}$ of $C^\infty$ diffeomorphisms $\Psi_x : \W(x) \to E_x := T_x\W(x)$ such that
\begin{enumerate}
    \item  $P_x = \Psi_{f(x)} \circ f \circ \Psi^{-1}_x: E_x \to E_{f(x)}$ is a polynomial map with degree at most $d$ for each $x\in M$;
    \item $\Psi_x(x) = 0$ and $D_x\Psi_x$ is the identity map for each $x\in M$,
    \item $\Psi_x$ depends continuously on $x\in M$ in the $C^\infty$ topology and depends smoothly on $x$ along the leaves of $\W$.
    \item For any $g\in \Z^0(f)$, if $g$ is uniformly $C^{s_0}$ along $\W$, then $\Psi_{g(x)} \circ g \circ \Psi^{-1}_x: E_x \to E_{g(x)}$ is also a polynomial map with degree at most $d$, for each $x\in M$. In particular, $g$ is $C^\infty$ when restricted to the leaves of $\W^{s/u}_f$ .
 
\end{enumerate}\end{lemma}
\begin{remark}\label{rema: dependence s0 f}
If $f$ has narrow  spectrum, i.e. $f$ has specified $(\chi,\epsilon)$-spectrum on $\W$, the parameter $s_0(f)$ can be taken to depend only on $\chi, \epsilon$. In particular, for $(\chi',\epsilon')$ sufficiently close to $(\chi,\epsilon)$, and for
any $f'$ that has $(\chi',\epsilon')$ spectrum on a contracting foliation $\W'$,  we have that  $s_0(f')$ is close to $s_0(f)$.  The same stability under $C^1$-small perturbations holds for the $1/2$-pinching condition. 
\end{remark}





\section{Proof of Theorem \ref{tmain=transitivecent}}
We return to the topic of diffeomorphisms with transitive centralizer.  Let $f\colon M\to M$ be a homeomorphism of a closed manifold $M$, and suppose that  $\Z^r(f)$ acts transitively on $M$, for some $r\geq 1$.
\subsection{A useful regularity result}
We have the following useful proposition about the regularity of mappings.

\begin{prop}\label{prop: reg prop}Let $r\geq 1$ be an integer, and let $g: M\to M$ be a continuous map on a compact manifold $M$. If for any $x,x'\in M$ there exists $\varphi\in \Z^r(g)$ such that $\varphi(x)=x'$, then $g$ is a $C^r$ map.
\end{prop}
\begin{proof}By continuity of $g$ and compactness of $M$, the graph $\Gamma_g:=\{(x,g(x)), x\in M\}$ of $g$ is a compact subset of $M\times M$. For any two points $(x,g(x)), (x',g(x'))$ in $\Gamma_g$, there exists $\varphi\in \Z^r(g)$ such that $\varphi(x)=x'$.  Then
\begin{eqnarray*}
\varphi\times\varphi(x,g(x)) &=& (\varphi(x), \varphi(g(x)))\\
&=& (x',g(\varphi(x))) = (x',g(x')).
\end{eqnarray*}
Since $\varphi\times\varphi$ is a $C^r$ diffeomorphism, this means that the compact set $\Gamma_g$ is $C^r$ {\it homogeneous} (in the sense of \cite[p.8]{W}). 
Every locally compact, $C^r$ homogeneous subset of a manifold is a $C^r$ submanifold \cite[Theorem B]{W} (see also \cite{RSS}); 
it follows that $\Gamma_g$ is a $C^r$ submanifold in $M\times M$, and hence the projection maps $\mathrm{Pr}_{1,2}|_{\Gamma_g}$ to each factor of $M\times M$ are $C^r$ when restricted to $\Gamma_g$. Moreover, $\mathrm{Pr}_1|_{\Gamma_g}$ is a $C^r$ homeomorphism. By Sard's theorem,  (see \cite{Hi}), there exists a point $(x,g(x))\in \Gamma_g$ such that $(x,g(x))$ is not a critical point of $\mathrm{Pr}_1|_{\Gamma_g}$. 

The inverse mapping theorem implies that $\mathrm{Pr}_1|_{\Gamma_g}$ has a $C^r$ inverse in a neighborhood of $x$, which implies that $g$ is $C^r$ in a neighborhood of $x$. Now for an arbitrary point $x'\in M$ and a sufficiently small neighborhood $B(x')\subset M$ of $x'$, we take a $C^r$ diffeomorphism $\psi\in \Z^r(g)$ such that $\psi(x)=x'$.  Then in a suitably small neighborhood $U$ of $x$,  we have 
\[g|_{B(x')}=\psi|_{g(x)}\circ g|_U \circ \psi^{-1}|_{x'},\]
which implies $g$ is also $C^r$ in a possibly smaller neighborhood of $x'$. As  $x'\in M$ was arbitrary, $g$ is $C^r$.
\end{proof}
\subsection{Nice subgroups and nice filtrations}
\begin{definition}\label{def: sgm cpct sbgp} Let $H$ be a topological group. A subgroup $\mathcal K\subset H$ is called \textit{nice} if there exists a countable increasing family of compact sets $K_i, i=1,2,\dots$ such that $\mathcal K=\cup_{i=1}^\infty K_i$, and for any $i,j\in \ZZ^+$, there exists $l=l(i,j)\in \ZZ^+$ such that $K_i\cdot K_j\subset K_l$ and $K_i^{-1}\subset K_l$. The family $\{K_i\}$ is called a \textit{nice filtration} of $\mathcal K$. It is easy to see a nice subgroup $\mathcal K$ may contain many nice filtrations.
\end{definition}

\noindent
\textbf{Examples:} 
\begin{itemize}
\item Let $M$ be a compact manifold and $H=\Ho(M)$ be the group of homeomorphisms of $M$. Then the subgroup of $H$ formed by the bi-Lipchitz homeomorphisms $\mathrm{Lip}(M)$ is nice since $\mathrm{Lip}(M)=\cup_{k=1}^\infty \mathrm{Lip}^k(M)$, where $\mathrm{Lip}^k(M)$ is the set of all bi-Lipchitz homeomorphism of $M$ with Lipchitz constant bounded by $k$.
\item Similarly, the subgroup of $H$ formed by all bi-H\"older homeomorphisms of $M$, $\mathrm{Hol}(M)=\cup_{k=1} \mathrm{Hol}^k(M)$ is nice, here $ \mathrm{Hol}^k(M)$ denote the set of bi-H\"older homeomorphisms with exponents bounded by $1/k$ and constants bounded by $k$.
\item  Lemma \ref{lemma: holo grp nice} gives examples arising in partially hyperbolic dynamics.  
\end{itemize}

Combining the nice property with Proposition \ref{prop: reg prop}, we have the following crucial proposition.

\begin{prop}\label{p=minimalhomogeneous} Let  $r\geq 0$ be an integer, and let $g\colon M\to M$ be a homeomorphism of a connected closed manifold $M$. Suppose that $\mathcal K$ is a nice subgroup of $\Z^r(g)$ such that for every $x,x'\in X$ there exists $\varphi\in \mathcal K$ such that $\varphi(x) = x'$.  Then 
the following holds:
\begin{enumerate}
\item $g$ is a $C^r$ diffeomorphism (a homeomorphism if $r=0$). 
\item \label{lem=continuityofvarphi} Let  $\{K_i\}_{i\in \ZZ^+}$ be a nice filtration  of $\mathcal K$.  There exists $N\geq 1$ such that for every $\varepsilon>0$, there is $\delta>0$ satisfying: for all $x,x'\in M$, $d(x,x')<\delta$ implies existence of $\varphi\in K_N$ with $d_{C^r}(\varphi,\id)<\varepsilon$ and  $\varphi(x)=x'$. 
\item \label{lem=equicontinuity} The family $\{g^n\colon n\in \ZZ\}$ is precompact in the $C^{0}$ topology. 
The action by $\G := \overline{\{ g^n : n\in \ZZ\}}^{C^{0}}$ is a free compact abelian group action on $M$. 
\item \label{Prop item upgr C0 to Cr}Moreover, if $g$ is minimal, then there exist a $C^r$ diffeomorphism $h: \TT^d\to M$ and $\rho\in \TT^d$ such that $h^{-1}\circ g\circ h(x) = x + \rho$, for all $x\in \TT^d$.
\end{enumerate}
\end{prop}

\noindent{\bf Question:} {\em Does the same result in (3) hold if we remove the $\sigma$-compact assumption and only assume $X$ is a compact metric space? That is, suppose that $g\colon X\to X$ is a minimal homeomorphism and that for every $x,x'\in X$, there exists a homeomorphism $h\colon X\to X$ such that $hg=gh$ and $h(x)=x'$.  Is $X$ a compact abelian group?}

\begin{proof}[Proof of Proposition~\ref{p=minimalhomogeneous}.] Since $\Z^r(g)=\Z^r(g^{-1})$,  Proposition \ref{prop: reg prop} and the transitivity of $\Z^r(g)$ imply that both $g$ and $g^{-1}$ are $C^r$, and so $g$ is a $C^r$ diffeomorphism.

Fix a nice filtration $\{K_i\}$ of $\mathcal K$.  Fix $z\in M$ and consider the evaluation map
$\Phi\colon  \Z^r(g) \to M$ defined by
\[\Phi(\varphi) =  \varphi(z).\]
The compactness of $K_i$ in $\Z^r(g)$ implies that for any $z\in M$, $K_i\cdot z$ is compact in $M$ (since the evaluation map is continuous)  A straightforward application of Baire category theorem gives the following.

\begin{lemma}\label{lem=nonemptyinterior}There exists $N_0>0$ such that for any $y,z\in M$, there exists $\varphi\in K_{N_0}$ such that $\varphi(z)=y$.
\end{lemma}

\begin{proof} Since $\mathcal K$ admits a nice filtration $\{K_i\}$, it suffices to show a slightly weaker version of this lemma:  for any $z\in M$, there exists $N_1\in \ZZ^+$ such that
for every $y\in M$, there exists $\varphi\in K_{N_1}$ such that $\varphi(z)=y$.

Fix $z\in M$. The transitivity of the $\mathcal K$ action implies that 
\[M = \bigcup_{i\in \NN} K_i\cdot z.
\]
is a countable union of closed subsets of a complete metric space $M$. The Baire category theorem implies that
there exists $i\in \NN$ such that $K_i\cdot z$ has nonempty interior in $M$.  Fix $z'\in M$ and $\delta>0$ such that $B(z',\delta) \subset K_i\cdot z$.  For $y\in M$, there exists $\varphi = \varphi_{z',y}\in \mathcal K$ such that $\varphi(z') = y$.  Then  there exists $N_{y}$ such that for every  $y'\in \varphi(B(z',\delta))$, there exists $\varphi'\in  K_{N_y}$ such that $\varphi'(z) = y'$: one just composes an element of $ K_i$ to go from $z$ to an element of $B(z',\delta)$ with $\varphi_{z',y}$.   One thus obtains an open cover of $M$;  extracting a finite subcover gives the result.
\end{proof}

Following the treatment in Avila-Santamaria-Viana \cite{ASV}, we consider a continuous map $\Phi\colon \A\to \B$  between topological spaces $\A$ and $\B$.  We say a point $x\in \A$ is {\em regular} if for every neighborhood $\V$ of $x$ we have $\Phi(x)\in\mathrm{int}\left(\Phi(\V) \right)$.  A point $y\in \B$ is a {\em regular value} of $\Phi$ if every point of $\Phi^{-1}(y)$ is regular.
We will use the following result from \cite{ASV}, which is a type of Sard's theorem for continuous maps.

\begin{prop}\label{prop=ASV}\cite[Proposition 7.6]{ASV} Let $\A$ be a compact metrizable space, and let $\B$ be a locally compact Hausdorff space.  If $\Phi\colon \A\to \B$ is continuous, then the set of regular values of $\Phi$ is residual. 
\end{prop}
Note that Proposition~\ref{prop=ASV} implies that for such a continuous map $\Phi$, either the image of $\Phi$ is meager or $\Phi$ has regular points.

\begin{lemma}\label{lem=existsregularpoint} Fix $z\in M$, and suppose that $K_i\cdot z$ has nonempty interior in $M$, for some $i\in \ZZ^+$.  Then there exists $\varphi_0 \in  K_i$ such that for every $C^r$ neighborhood $\V$ of $\varphi_0$ in  $K_i$, we have
\[\varphi_0(z) \in \mathrm{int}\left(\V\cdot z \right).\]
\end{lemma}

\begin{proof} We consider the evaluation map
$\Phi\colon   K_i \to M$ defined by
\[\Phi(\varphi) =  \varphi(z).\]

It is a continuous map from a compact metric space to a compact Hausdorff space.  Proposition~\ref{prop=ASV} implies that the set of regular values of $\Phi$ is residual.   Let $y_0\in K_i$ be a regular value, and let $\varphi_0\in \Phi^{-1}(y_0)$ be a regular point.  Then by the definition of regular point, $\varphi_0$ satisfies the desired property.
\end{proof}

Now we show \eqref{lem=continuityofvarphi} in Proposition \ref{p=minimalhomogeneous}. Fix $z\in M$, and let $N_0\in \ZZ^+$ is given by Lemma \ref{lem=nonemptyinterior}, then by Lemma~\ref{lem=existsregularpoint}, we can pick $\varphi_0\in K_{N_0}$ and $y_0:=\varphi_0(z)$ satisfying Lemma~\ref{lem=existsregularpoint}.  Given $\varepsilon>0$, let  $\V$ be a neighborhood of $\varphi_0$ in $K_{N_0}$ such that
\[d_{C^r}(\V^{-1} \V,\id)<\varepsilon.
\]
Note that since $\varphi_0$ satisfies the conclusions of Lemma~\ref{lem=existsregularpoint}, we have that  $\V \cdot z$ contains a neighborhood $W_{y_0}$ of $y_0$.

Given $y\in M$, choose $\varphi\in K_{N_0}$ such that $\varphi(y_0)=y$.  Then $\varphi(\V \cdot z)$ contains a neighborhood $W_y=\varphi(W_{y_0})$ of $y$.
We thus obtain an open cover $\{W_y : y\in M\}$ of $M$; let $\delta>0$ be the Lebesgue number of this cover.

Now given $x,x'\in M$, if $d(x,x')<\delta$, there exists a $W_y$ such that $x,x'\in W_y$.  This means that there exist $\varphi, \varphi'\in \V$ such that
$\varphi\varphi_0(z) = x$ and $\varphi'\varphi_0(z) = x'$.  In other words, $x' =  (\varphi')^{-1}\varphi(x)$.  Since  $(\varphi')^{-1}\varphi\in \V^{-1}\V$, we have that $d_{C^r}( (\varphi')^{-1}\varphi,\id) < \varepsilon$.   This completes the proof of \eqref{lem=continuityofvarphi}.

Now we show \eqref{lem=equicontinuity} in Proposition \ref{p=minimalhomogeneous}. Given $\epsilon>0$, let $\delta>0$ be given by \eqref{lem=continuityofvarphi}.  If $d(x,x')<\delta$, then there exists $\varphi\in K_{N_0}$ with $d_{C^r}(\varphi^{\pm 1},\id)<\varepsilon$ such that $\varphi(x)=x'$.  

Then for any $n\in \ZZ$, \begin{eqnarray*}d(g^n(x), g^n(x')) &=& d(g^n(x), \varphi \circ g^n\circ \varphi^{-1}(x')) \\ &<&\varepsilon, 
\end{eqnarray*}
since both  $\varphi^{-1}$ and $\varphi$ are $C^r$-close to the identity in a neighborhood of $x'$ and $x$, respectively.
The Arzel\`a-Ascoli theorem then implies that  $\{g^n:n\in \ZZ\}$ is precompact in $C^{0}(M)$. 
 Let $\G$ be the $C^{0}$ closure of the iterates of $g$:
\[\G := \overline{\{ g^n : n\in \ZZ\}}^{C^{0}} \subset \Ho(M).
\]
By definition of $\G$ it clearly commutes with any element of $\mathcal Z^r(g)$, therefore $\G$ commutes with $\mathcal K$. Since $\mathcal K$ acts transitively on $M$ and commutes with $\G$, if $g'\in \G$ fixes a point $x\in M$, then $g'$ fixes all points and is the identity. Therefore the $\G$ action is free. This 
complets the proof of \eqref{lem=equicontinuity}.

To show part \eqref{Prop item upgr C0 to Cr} of Proposition \ref{p=minimalhomogeneous}, first we consider its $C^0$ version. Assume that $g$ is minimal. Using the $C^0$  precompactness  of $\{g^n:n\in \ZZ\}$ from \eqref{lem=equicontinuity}, we construct an abelian topological group structure on $M$.  As before, let 
$\G := \overline{\{ g^n : n\in \ZZ\}}^{C^{0}}$.


Minimality of $g$ implies that  $\G$ acts transitively on $M$: given $x,y\in M$, since $g$ is minimal, 
there exists  $n_j\to \infty$  such that $g^{n_j}(x)\to y$; using compactness of $\G$ and passing to a subsequence we obtain in the limit a map $\bar{g} = \lim_k g^{n_{j_k}} \in \G$ satisfying $\bar{g}(x) = y$.

Fix $z\in M$. Since the $\G$-action is transitive and free, the evaluation map $\Phi\colon \G\to M$ defined by $\Phi(\varphi) = \varphi(z)$ is a homeomorphism.  It defines a continuous abelian group structure on $M$ via the operation $\varphi_1(z) + \varphi_2(z) := \varphi_1(\varphi_2(z)) = \varphi_2(\varphi_1(z))$.  Thus $M$ is topologically a torus, and $\varphi\in \G$ is topologically conjugate via $\Phi$ to a minimal translation $x\mapsto x+\rho$ for some $\rho\in \TT^d$.



Proposition \ref{prop: reg prop} implies that the elements of $\G$ are $C^r$.  Then \cite[Theorem 5]{CM} implies that the action by $\G$ is in fact  a $C^r$  Lie group  action, and $M$ is $C^r$ conjugate to the standard torus. 

\end{proof}
\subsection{Proof of Theorem \ref{tmain=transitivecent}}
\begin{proof}[Proof of Theorem \ref{tmain=transitivecent}] (1) First we assume $\Z^r(f)$ acts transitively on $M$. Proposition \ref{prop: reg prop} implies that $f$ is $C^r$. Let $\mathcal K$ be the group of all  bi-Lipschitz homeomorphisms of $M$ that commute with $f$, and note that $\mathcal K$ is a nice subgroup of $\Ho(M)$ (since any $C^0$ limit point of a sequence of bi-Lipchitz homeomorphisms with Lipchitz constant bounded by $L$ and commuting with $f$, is bi-Lipchitz and with the same Lipchitz constant and commutes with $f$). Then by Proposition \ref{p=minimalhomogeneous}, $\{f^n,n\in \ZZ\}$ is precompact in $\Ho(M)$, and $P=\overline{\{f^n,n\in \ZZ\}}$ is a compact topological group acting on $M$. Since  each $p\in P$  commutes with $\Z^r(f)$ which acts transitively on $M$,  Proposition \ref{prop: reg prop} implies that $p$ is $C^r$. Then $P$ is a compact topological group acting continuously on $M$ by $C^r$ diffeomorphisms, and so \cite[Theorem 5]{CM} implies that $P$ is a Lie group and its action is a $C^r$ compact Lie group action on $M$. 

Since $P$ is the $C^0$-closure of $\{f^n, n\in \ZZ\}$, it follows that $P$ is a compact abelian subgroup of $\Ho(M)$, containing a dense subset $\{f^n, n\in \ZZ\}$. Therefore $P$ is a direct product of some $\ZZ/d\ZZ$ with $\TT^k$. 

Since $P$ is compact Lie group, its action is proper. By part \eqref{lem=equicontinuity} of Proposition \ref{p=minimalhomogeneous}, the $P$ action is free.  
Then we apply the classical fact that a proper and free Lie group action induces a $C^r$ principal bundle structure on $M$. This implies the ``only if" part of Theorem \ref{tmain=transitivecent}.


(2) Now we assume that $M$ is a $C^r$ principal $\TT^k\times \ZZ/d\ZZ$ bundle and $f$ is a constant minimal translation on the bundle and identity on the base.  Let 
$G\to M \xrightarrow{\pi}\overline{M}$ be the $f$-invariant principal fiber bundle structure, with $G= \mathbb T^k\times \ZZ/d\ZZ$.  Since $M$ is connected, to show that $\Z^r(f)$ acts
 transitively on $M$, it  suffices  to show for any $x\in M$ and any $y$ close to $x$, there exists $g\in \Z^r(f)$ such that $g(x)=y$. To show this we take a small neighborhood $\bar{U}$ of $\pi(x)$ such that $\pi^{-1}(\bar U)$ is $C^r$ diffeomorphic to $\bar{U}\times E$. Then for any $y\in \pi^{-1}(\bar U)$, it is not hard to construct a $C^r$ diffeomorphism $g$ of $M$ such that
 \begin{enumerate}
     \item $g$ is supported on $\pi^{-1}(\bar U)$.
     \item $g$ preserves the principal bundle structure in $\pi^{-1}(\bar U)$, i.e. $g$ is translation of $G$ along each fiber.
     \item $g(x)=y$ (such a $g$ is constructed as a translation skew product of a map $\bar{g}:\overline{M}\to \overline{M}$ supported on $\bar{U}'$ such that $\overline{\bar{U}'}\subset \bar{U}$).
     \end{enumerate}
Since $f$ is fiber fixing and a constant translation along every fiber, it follows that $fg=gf$.
\end{proof}

\subsection{A byproduct of the proof of Theorem \ref{tmain=transitivecent}}
For later use we consider the following corollary of Propositions~\ref{prop: reg prop} and \ref{p=minimalhomogeneous}. For a set $\mathcal K\subset \Ho(M)$, and $k\in \NN$, we define \[\Z^k(\mathcal K):= \{g\in \Diff^k(M), gg'=g'g,~~~ \forall g'\in \mathcal K\}.\] 

\begin{maintheorem}\label{mainthm: recstrct Lie grp actn} Let $M$ be a connected closed manifold and let $\mathcal K$ be a nice subgroup of $\Ho(M)$ such that $\mathcal K\subset \Diff^r(M)$, for some $r\geq 1$. If  $\mathcal K$ acts transitively on $M$, then $M$ is a $C^r$ principal fiber bundle (possibly a trivial one). Moreover $\Z^0(\mathcal K)=\Z^r(\mathcal K)$ consists of the set of constant translations along the fibers (that is, the right action of the structure group).
\end{maintheorem}
\begin{proof}Using the same method as in the proof of \eqref{lem=equicontinuity} of Proposition \ref{p=minimalhomogeneous},  we use the Arzel\`a-Ascoli theorem to  conclude that $\Z^0(\mathcal K)$ is a compact group (in the $C^0$ topology). Transitivity of the action of $\mathcal K$ implies that $\Z^0(\mathcal K)=\Z^r(\mathcal K)$. This implies that $\Z^0(\mathcal K)$ is a Lie group\cite{RS}. Then \cite[Theorem 5]{CM} implies that the $\Z(\mathcal K)$ action is a $C^r$ action of a Lie group on $M$. 

For any $g\in \Z^0(\mathcal K)$, if $g$ has a fixed point $x$, by transitivity of $\mathcal K$ action we can easily show that any point $x'$ of $M$ is also a fixed point of $g$, which means that $g$ is the identity. Thus the action of $\Z^0(\mathcal K)$ is a free and proper hence induces a $C^r$ principal $\Z^0(\mathcal K)$-bundle structure on $M$.
\end{proof}

 \section{Applications to fibered partially hyperbolic systems}
 \subsection{The center fixing centralizer}
In the following theorem, the integer parameter $s_0\geq 1$ controls the spectral behavior of $Df|_{E^{s,u}}$, and $r$ measures in part the nonconformality of $Df|_{E^c}$ (the larger the $r$, the more conformal). So in general $r$ can be arbitrarily large, i.e. close to or equal $\infty$, for example, when $f$ is a perturbation of an isometric extension. In particular, $r$ could be much larger than $s_0$, and in that case the following theorem demonstrates a certain bootstrapping effect on regularity.
\begin{maintheorem}\label{t=torusfibervp}  Let $f\in \Diff^\infty(M)$ be a fibered partially hyperbolic diffeomorphism that is accessible and $r$-bunched, for some $r\in \ZZ^+$. We further assume that the cocycles $Df|_{E^s}, Df^{-1}|_{E^u}$ either satisfy the the pointwise $1/2$-pinching condition or have narrow spectrum. Then there exists $s_0\geq 1$ such that for any $s\in \{s_0,\ldots, r\}$,  $\CZ^s(f)$ is a compact $k$-dimensional Lie subgroup of $\Diff^s(M)$, for some $k\leq \dim E^c_{f_0}$, such that 
\begin{enumerate}    
    \item The action by $\CZ^s(f)$ is a $C^r$ Lie group action that is uniformly $C^\infty$ along the stable and unstable foliations of $f$.
    \item $M$ admits an $f$-invariant, holonomy-invariant $C^r$-principal bundle structure that is subordinate to $\W^c_f$.
    \item Moreover, if $\dim \CZ^s(f)=\dim \W^c_f$, then  $\CZ^s(f)=\CZ^\infty(f)$ and $f$ is a $C^\infty$ compact Lie group extension of an Anosov diffeomorphism. 
\end{enumerate}   
\end{maintheorem}




\begin{remark}
From the proof we will see our result also holds for $f\in \Diff^k(M)$ for some $k<\infty$ sufficiently large; for simplicity we only state the $C^\infty$ case. We have the following corollary of the proof of Theorem~\ref{t=torusfibervp}.
\end{remark}

\begin{coro}\label{coro: 1d center}
 Let $f\in \Diff^\infty(M)$ be an accessible fibered partially hyperbolic diffeomorphism with $1$-dimensional center. We further assume that the cocycles $Df|_{E^s}, Df^{-1}|_{E^u}$ either satisfy the the pointwise $1/2$-pinching condition or have narrow spectrum. Then $\CZ^\infty(f)$ is either finite or $\TT^1$. In the latter case $M$ admits an $f$-invariant, holonomy-invariant $C^\infty$ principal bundle structure along $\W^c$. 
\end{coro}

 We will prove  Theorem~\ref{t=torusfibervp} and Corollary \ref{coro: 1d center}  in the following subsection; we first establish some preliminaries. 
 
\begin{lemma}\label{lemma: holo grp nice}
Let $r\in \ZZ^+$ and $f$ be $C^{r}$ a fibered partially hyperbolic diffeomorphism that is accessible and $r$-bunched. Then for any fixed center leaf $\W^c_f(x_0)$, the set of all $su$-holonomies from $\W^c(x_0)$ to $\W^c(x_0)$ is a nice subgroup of $\Diff^k(\W^c(x_0))$ for $k=0,\dots, r$.
\end{lemma}
\begin{proof}[Proof of Lemma~\ref{lemma: holo grp nice}] If $f$ is  $C^{r}$ and $r$-bunched, then  the stable and unstable holonomies between center leaves are $C^r$ and vary continuously  in $C^r$ the topology as a function of the paths inducing them.  Therefore the set of $su$-holonomies from a fixed center leaf $\W^c(x_0)$ to itself admit a nice filtration by an increasing sequence of compact sets $K_n\subset \Diff^k(\W^c(x_0))$, where $K_n$ is the subset of $\Diff^r(\W^c(x_0))$ induced by  $su$-paths  with  $\leq n$ legs each of length $\leq n$.
\end{proof}

It is easy to see that the holonomy maps of a fibered partially hyperbolic diffeomorphism $f\in \Diff^1(M)$ commute with the elements of $\CZ^1(f)$.  This is not necessarily the case for elements of $\CZ^0(f)$, as homeomorphisms commuting with $f$ might not preserve the leaves of $\W^{s/u}_f$.  This motivates the following definition.

\begin{definition}\label{def: CZ_0 redef}
For a fibered partially hyperbolic diffeomorphism $f\in \Diff^1(M)$, we denote by 
\begin{enumerate}
    \item $\CZ_\ast(f)$ the group of homeomorphisms of $M$ that commute with $f$, fixing each $\W^c_f$ leaf and preserving the $\W^{s/u}_f$-foliations; and
    \item $\C^r\Z_\ast(f)\subset \CZ_\ast$ the subgroup of $\CZ_\ast(f)$ consisting of the elements that are uniformly $C^r$ on $\W^c_f$ leaves. 
\end{enumerate} 
Clearly we have $\CZ^r(f)\subset \C^r\Z_\ast(f) \subset \CZ_\ast(f)$.
\end{definition}

We have the following proposition, which is an application of Theorem \ref{mainthm: recstrct Lie grp actn}.
\begin{prop}\label{prop: coro Tm 10} Let  $f\in \Diff^\infty(M)$ be a fibered partially hyperbolic diffeomorphism. If $f$ is accessible and $r$-bunched, for some $r\in \ZZ^+$, then
\begin{enumerate}
    \item $\C^r\Z_\ast(f) =\CZ_\ast(f)$, i.e., every element in $\CZ_\ast(f)$ is uniformly $C^r$ along $\W^c_f$.
    \item The action of $\CZ_\ast(f)$ on $M$ is a free compact Lie group  action. 
\end{enumerate} 
\end{prop}
\begin{proof}[Proof of Proposition \ref{prop: coro Tm 10}]
 For any fixed center leaf $\W^c_f(x_0)$, consider the group $\mathcal K(x_0)$ formed by $su$-holonomies from $\W^c_f(x_0)$ to $\W^c_f(x_0)$. Lemma~\ref{lemma: holo grp nice} implies that $\mathcal K(x_0)\subset \Diff^r(\W^c_f(x_0))$ is a nice subgroup of $\Ho(\W^c_f(x_0))$.  Theorem~\ref{mainthm: recstrct Lie grp actn} then implies that any element in $\CZ_\ast(f)$ is $C^r$ along $\W^c_f(x_0)$.
 
Claim (1) of Proposition \ref{prop: coro Tm 10} is implied by the following \textit{uniformity lemma}. 

\begin{lemma}\label{lemma: unfmty lm} 
Let $f$ be $C^{r}$ fibered partially hyperbolic diffeomorphism that is  $r$-bunched, for some $r\in \ZZ_+$. Suppose that $g\in \CZ_\ast(f)$ is $C^r$ along a center leaf 
 $\W^c_f(x_0)$. Then $g$ is uniformly $C^r$ along all leaves of $\W^c_f$. Moreover, $y\mapsto g\vert_{\W^c_f(y)}$ induces a continuous map from $M$ into the space of $r$-jets along leaves of $\W^c_f$.
\end{lemma}
\begin{proof}[Proof of Lemma~\ref{lemma: unfmty lm}]  
 For $x,y\in M$, let $\gamma$ be an $su$-path from $x$ to $y$. The induced holonomy $H_\gamma\colon \W^c_f(x)\to \W^c_f(y)$ is $C^r$ (since $f$ is assumed to be $r$-bunched), and it varies continuously in $C^r$ topology with the path $\gamma$. For $g\in \CZ_\ast(f)$, we have that $H_\gamma\circ g\vert_{\W^c_f(x)}\circ H_\gamma^{-1} = g\vert_{\W^c_f(y)}$. 

Let $d=\dim M-\dim E^c$. Due to the fibered structure of $\W^c_f$, for any $x\in M$, there is a continuous family of $2$-legged $su$-paths  $\{\gamma_\xi\,\mid\,\xi\in [-\epsilon, \epsilon]^{d}\}$ starting at $x$ such that 
the map $\W^c_f(x)\times [-\epsilon, \epsilon]^{d}\to M$ defined by $(\xi, z)\mapsto H_{\gamma_{\xi}}(z)$ is a homeomorphism onto a  neighborhood  of $\W^c(x)$. This  defines a continuous map from  $[-\epsilon, \epsilon]^{d}$ into the space of $C^r$ embeddings of $\W^c_f(x)$ into $M$ via $\xi\mapsto H_{\gamma_{\xi}}$.

If $g\vert_{\W^c_f(x)}$ is $C^r$, then $y\mapsto g\vert_{\W^c_f(y)}$  is continuous in the $C^r$ topology in this neighborhood of $g\vert_{\W^c_f(x)}$, via the parametrization $\xi\mapsto H_{\gamma_\xi}\circ g\vert_{\W^c_f(x)}\circ H_{\gamma_\xi}^{-1}$. Since $M$ is connected, if $g\vert_{\W^c_f(x_0)}$  is $C^r$ for some $x_0$, then  $y\mapsto g\vert_{\W^c_f(y)}$  is continuous in the $C^r$ topology on all of $M$.
\end{proof}

Now we prove part (2) of Proposition \ref{prop: coro Tm 10}. We first show that the action  of $\CZ_\ast$ on $M$ is free. To this end, suppose that
 $g\in \CZ_\ast$ fixes a point $x$. Then for every point $y\in M$, there is an $su$ path $\gamma$ from $x$ to $y$, and $g(y) = H_\gamma \circ g \circ H_\gamma^{-1}(y) = y$.  Thus $g$ is the identity, and the action is free.
 
The proof of Lemma~\ref{lemma: unfmty lm}, shows that for any $x_0\in M$, the group representation $\iota: \CZ_\ast\to \Ho(\W^c_f(x_0))$ defined by restriction to $\W^c_f(x_0)$ is a topological group embedding. Lemma~\ref{lemma: holo grp nice} implies that $\mathcal K(x_0)$ is a nice subgroup of $\Ho(\W^c_f(x_0))$,  which centralizes $\iota(\CZ_\ast)$.
Accessibility of $f$ implies that $\mathcal K(x_0)$ acts transitively on $\W^c_f(x_0)$. Theorem \ref{mainthm: recstrct Lie grp actn} then implies that $\iota(\CZ_\ast)$ is contained in the compact Lie group 
$\Z(\mathcal K(x_0)):=\{g\in \Ho(\W^c_f(x_0)): gh=hg,~~ \forall h\in \mathcal K(x_0)\}$. 

Since $\iota$ is faithful, to show the action by $\CZ_\ast$ is a compact Lie group action,  we need only show that $\iota(\CZ_\ast)$ is closed in $\Ho(\W^c_f(x_0))$. If not, there exists a sequence $g_n\in \CZ_\ast$ such that $g_n$ has no limit point in $\Ho(M)$ but the restriction $g_n|_{\W^c_f(x_0)}$ converges uniformly. But the proof of Lemma~\ref{lemma: unfmty lm} shows that uniform convergence of $g_n|_{\W^c_f(x_0)}$ implies uniform convergence of $g_n$ on $M$, a contradiction. This establishes part (2) of of Proposition \ref{prop: coro Tm 10}, completing its proof.
\end{proof}

\subsection{Proof of Theorem~\ref{t=torusfibervp} }
Assume that $f$ satisfies the hypotheses of Theorem \ref{t=torusfibervp}.  We denote $\CZ^s(f)$ by $\CZ^s$, and let $\CZ_\ast = \CZ_\ast(f)$  and $\CZ^r_\ast=\CZ^r_\ast(f)$ be defined as in the previous subsection.  Since $Df|_{E^s}$ and $Df^{-1}|_{E^u}$ satisfy either the $\frac{1}{2}$-pinching or narrow spectrum condition, part (4) of Lemma \ref{lemma: normal form} implies that there exists $s_0>0$,  such that if $s\geq s_0$, then any $g\in \CZ^s$ is $C^\infty$ along the leaves of $\W^{s/u}$ (and is polynomial after a uniformly $C^\infty$ change of coordinates along leaves).  

Moreover,  Remark \ref{rema: dependence s0 f} (see \cite{Kal-normalform} for more details), implies that  by increasing $s_0(f_0)$ if necessary, we can take $s_0(f_0)=s_0(f)$ to work for all $f$ such that $d_{C^1}(f_0,f)<\varepsilon$, for some small $\varepsilon$.  We fix such an $s_0$.  Proposition~\ref{prop: coro Tm 10} gives that $\CZ^s\subset \CZ_\ast=\C^r\Z_\ast$ for $s_0\leq  s\leq r$, and in particular, every element of $\CZ^s$ is $C^r$ along the leaves of $\W^{c}_f$. Since, in addition,  the elements of $\CZ^s$ are $C^\infty$ along the $\W^{s/u}_f$  leaves,  Journ\'e's lemma \cite{Journe} implies that $\CZ^s=\CZ^r$.  Thus any element of $\CZ^s$ is $C^\infty$ along the leaves of $\W^{s/u}_f$.

To show the action by $\CZ^s$ is a $C^r$ compact Lie group action, we recall from  Proposition \ref{prop: coro Tm 10} that the action by $\CZ_\ast$ is a free compact Lie group action. It then suffices to show the following.
\\
\\
\textbf{Claim:}
The group $\CZ^s$ is a $C^0$ closed subgroup of $\CZ_\ast$. 
\\
\\
The claim implies that $\CZ^s$ is a compact Lie subgroup of $\CZ_\ast\subset \Ho(M)$, and each element acts on $M$ by $C^r$ diffeomorphisms; it follows from \cite[Theorem 5]{CM} that $\CZ^s$ acts by a $C^r$ Lie group action.
\begin{proof}[Proof of claim]  Consider an arbitrary sequence $g_n\in \CZ^s$ converging uniformly to $g_0\in \CZ_\ast$.  Proposition~\ref{prop: coro Tm 10} implies that $g_0\in \CZ^r_\ast$, and so $g_0$ is uniformly $C^r$ along $\W^c_f$-leaves.     

Part (4) of Lemma \ref{lemma: normal form} implies that, up to a uniformly $C^\infty$ coordinate change along $\W^{s/u}$, each $g_n$ is a polynomial map along $\W^{s/u}$ with bounded degree and bounded coefficients;
a $C^0$ limit of polynomial maps with bounded degree and bounded coefficients is polynomial, and so therefore, under a uniform $C^\infty$ coordinate change, the map $g_0$ is also a polynomial map along $\W^{s/u}$. Thus $g_0$ is uniformly $C^\infty$ along $\W^{s/u}$. Then by Journ\'e's lemma \cite{Journe}, $g_0$ is a $C^r$ diffeomorphism, and so $g_0\in \CZ^s$.
\end{proof}



 As in the proof of Theorem~\ref{mainthm: recstrct Lie grp actn}, we obtain from the action of $\CZ^s$ a principal bundle structure on $M$. 
  It is $f$-invariant, holonomy invariant and subfoliates $\W^c_f$; this gives (2) of Theorem \ref{t=torusfibervp}. Clearly $\dim \CZ^s$ is well-defined and not greater than $\dim W^c_f$ (a Lie group cannot act freely on a manifold with lower dimension).

To complete the proof of Theorem \ref{t=torusfibervp}, we consider the case that $\dim \CZ^s=\dim \W^c_f$. Denote by  $cc_{\id}(\cdot)$ the connected component of the identity in a given topological group.   
Since the action by $\CZ_\ast$ is free, so is the $\CZ^s$ action. Therefore the action of $\CZ^s$ is a free $C^r$ compact Lie group action.
  By compactness and connectedness of $\W^c_f$ and  $cc_{\id}(\CZ^s)$,  the group  $cc_{\id}(\CZ^s)$ acts freely and transitively on each $\W^c_f$ leaf, and by freeness of the action of $\CZ^s$, we have $\CZ^s=cc_{\id}(\CZ^s)$. In this case $\W^c_f$ is subfoliated by a principal $\CZ^s$-bundle, where $\CZ^s$ has the same dimension as $\W^c_f$. 

 Consider the bi-invariant metric on $cc_{\id}(\CZ^s)$: it  induces an $f$-invariant metric on $\W^c_f$. Then $f$ is $\infty$-bunched. But this implies that 
$\CZ^s=\CZ^\infty$, and $\W^c_f$ leaves are the orbits of the $C^\infty$, free compact Lie group action by $\CZ^s$.
To summarize, in the case that $\dim \CZ^s=\dim \W^c_f$,  the connected compact Lie group $\CZ^s$ acts on $M$ smoothly, freely and properly with orbit foliation $\W^c_f$. Thus the projection $\pi: M\to M/\W^c_f$ gives a $C^\infty$ principal bundle with structure group $\CZ^s$.   In any group the centralizer of all left translations is the group of right translations. Since $f$ commutes with $\CZ^s$,  it follows that
$f$ is a $C^\infty$ compact Lie group extension of an Anosov diffeomorphism. This completes the proof of Theorem~\ref{t=torusfibervp}.\qed

Using almost the same method, we can prove Corollary \ref{coro: 1d center}. For completeness  we give a sketch of the proof.  Proposition~\ref{prop: coro Tm 10} and Corollary \ref{coro: 1d center} imply that $\C^1\Z_\ast(f)$ is either finite or $\TT^1$. If it finite, then $\CZ^\infty(f)$ is obviously finite. If $\C^1\Z_\ast(f)$ is $\TT^1$, then $f$ is $\infty$-bunched. Then by almost the same proof as of Theorem \ref{t=torusfibervp}, using Journ\'e's lemma and normal form theory we obtain  that $\C^1\Z_\ast(f)=\C^\infty\Z_\ast(f)=\CZ^\infty(f)$ is $\TT^1$. The rest of Corollary \ref{coro: 1d center} is a consequence of part (3) of Theorem \ref{t=torusfibervp}.


\subsection{Proof of Theorem~\ref{main: application iso ext}}

Let  $f_0$ be  $\infty-$bunched, and assume that $Df_0|_{E^s}$, $Df_0^{-1}|_{E^u}$ have $1/2$-pinched or 
narrow spectrum. Lemma~\ref{lemma: normal form} gives an $s_0(f_0)\geq 1$,  such that for any $s\geq s_0$, the maps $f_0$  and any  $g\in \Z^{s}(f_0)$  preserve  $C^\infty$ normal forms along $\W^{s/u}_{f_0}$, and the same holds for any $f$ sufficiently $C^1$ close to $f_0$.  Fix $s\geq s_0$.  If  $f$ is sufficiently $C^1$-close to $f_0$, then $f$  is $r-$bunched for some $r>s$, is fibered (by \cite{HPS}) and accessible (by Lemma \ref{prop: stab acc}). Therefore such an $f$ satisfies all the assumptions of Theorem \ref{t=torusfibervp}.

Theorem \ref{t=torusfibervp} implies that $\CZ^{s}(f)$ is a $C^0$-closed Lie subgroup of $\Diff^s(M)$ with dimension $k \leq \dim E^c_{f_0}$.
Thus to show that $\Z^s(f)$ is also a $C^0$-closed Lie subgroup of $\Diff^s(M)$, it suffices to show that $\Z^s(f)$ is just countably many discrete copies of $\CZ^s(f)$ in $\mathrm{Homeo}(M)$. 

 Proposition \ref{prop: fib systm zf inv} implies that any element in $\Z^s(f)$ preserves the fiber bundle structure, and so  $\CZ^s(f)$ is a normal subgroup of $\Z^s(f)$, and the quotient group $\Z^s(f)/\CZ^s(f)$ has a continuous (with respect to the quotient metric) identification with a subgroup of $\Z^0(\bar{f})$, through the map $\iota: g\cdot \CZ^s(f)\mapsto \bar{g}$. 
 
The first paragraph of the proof of Proposition \ref{prop: fib systm zf inv} gives that $\bar{f}$ is H\"older conjugate to an Anosov automorphism on a nilmanifold $M$. Any transitive Anosov diffeomorphism (in particular an affine Anosov map on a compact nilmanifold) has countable discrete $C^0$ centralizer. Hence the coset partition $\Z^s(f)=\bigsqcup g\;\CZ^s(f)$ is actually formed by countably many discrete copies of $\CZ^s(f)$. Then $\Z^s(f)$ is a closed Lie group in $\Ho(M)$ and we have $\dim \Z^s(f)=\dim \CZ^s(f)$. The rest of results in Theorem \ref{main: application iso ext} are direct consequences of the corresponding results in Theorem \ref{t=torusfibervp} (for the spectral property of the base dynamics in (2), we use the fact that if $Df|_{E^s}$ has $(\chi,\epsilon)$-spectrum, and $\W^c_f$ is smooth then $D\bar{f}|_{E^s_{\bar{f}}}$ has $(\chi,\epsilon)$-spectrum as well, and similarly for the $1/2$-pinched case). This completes the proof of Theorem~\ref{main: application iso ext}.\qed
\section{Proof of Theorem \ref{coro: global Heis3}}
In this section we prove  Theorem \ref{coro: global Heis3}. Let $f$ be a  partially hyperbolic diffeomorphism on a $3$-dimensional non-toral nilmanifold $M$. Then $f$ satisfies the following properties:
\begin{enumerate}
    \item By \cite[Theorem 1.6]{HP}, $f$ is a fibered partially hyperbolic system over an Anosov homeomorphism $\bar f:N/\W^c_f\to N/\W^c_f$, which is  topologically conjugate to a hyperbolic automorphism of $\TT^2$.
\item By \cite[Proposition 6.4]{HP}, $f$ is accessible.
\item Since $E^{s/u}$ are $1$-dimensional, the restrictions $Df|_{E^{s/u}}$ satisfy the pointwise $1/2$-pinching condition. 
\end{enumerate}
Corollary \ref{coro: 1d center} then implies that $\CZ^{\infty}(f)$ is either finite or virtually $\TT^1$. In the latter case, $f$ preserves a $C^\infty$ principal $\TT^1$ bundle structure. 

To complete the proof of Theorem \ref{coro: global Heis3}, it suffices to show that $\Z^\infty(f)$ 
is virtually $\{f^n\circ \CZ^\infty(f), n\in \ZZ\}$.
Notice that $\Z^\infty(f)/\CZ^\infty(f)$ can be naturally embedded into $\Z^0(\bar{f})$. Since $\bar{f}$ is topologically conjugate to a toral hyperbolic automorphism, \cite{AP} and \cite[Proposition 3.7]{KKS} imply that $\Z^0(\bar{f})$ is virtually trivial, which forces $\Z^\infty(f)$ to be virtually $\{f^n\circ \CZ^\infty(f), n\in \ZZ\} \cong \ZZ\times \TT$.
\color{black}

\end{document}